\documentclass[10pt,twoside]{amsart}
\usepackage{geometry}
\usepackage[english]{babel}
\usepackage{graphicx}
\usepackage{amsmath}
\usepackage{amsfonts}
\usepackage{amssymb}

\newcommand{\wk}{\mbox{$\,<$\hspace{-5pt}\footnotesize )$\,$}}

\begin{document}

\numberwithin{equation}{section}
\newtheorem{teo}{Theorem}
\newtheorem{lemma}{Lemma}
\newtheorem{defi}{Definition}
\newtheorem{coro}{Corollary}
\newtheorem{prop}{Proposition}
\newtheorem{remark}{Remark}
\newtheorem{scho}{Scholium}
\numberwithin{lemma}{section}
\numberwithin{prop}{section}
\numberwithin{teo}{section}
\numberwithin{defi}{section}
\numberwithin{coro}{section}
\numberwithin{figure}{section}
\numberwithin{remark}{section}
\numberwithin{scho}{section}

\title{Geometric Properties of a Sine Function extendable to arbitrary normed planes}

\author{Vitor Balestro, Horst Martini, and Ralph Teixeira}
\address [V. Balestro] {CEFET/RJ Campus Nova Friburgo - Nova Friburgo - Brazil
\& Instituto de Matem\'{a}tica e Estat\'{i}stica - UFF - Niter\'{o}i - Brazil}
\email{vitorbalestro@mat.uff.br}
\address [H. Martini] {Fakult\"at f\"ur Mathematik - Technische Universit\"at Chemnitz - 09107 Chemnitz - Germany}
\email{martini@mathematik.tu-chemnitz.de}
\address [R. Teixeira] {Instituto de Matem\'{a}tica e Estat\'{i}stica - UFF - Niter\'{o}i - Brazil}
\email{ralph@mat.uff.br}

\begin{abstract} In this paper we study a metric generalization of the sine function which can be extended to arbitrary normed planes. We derive its main properties and give also some characterizations of Radon planes.
 Furthermore, we prove that the existence of an angular measure which is ``well-behaving'' with respect to the sine is only possible in the Euclidean plane, and we also define some new constants that estimate how non-Radon or non-Euclidean a normed plane can be. Sine preserving self-mappings are studied, and a complete description of the linear ones is given. In the last section we exhibit a version of the Law of Sines for Radon planes.
\end{abstract}

\maketitle

\section{Introduction} \label{secintro}

The notion of angle plays an important role in the geometry of (Minkowski or) normed planes. In \cite{brass} Brass introduced the notion of angular measure in normed planes and used such a measure to study packings of the unit circle. Later this concept was more explored. In \cite{duvelmeyer1} D\"{u}velmeyer characterized the Euclidean plane as the only plane where Busemann and Glogovskii angular bisectors can be defined by using an angular measure. Fankh\"{a}nel (see \cite{fankhanel1} and \cite{fankhanel2}) studied particular types of angular measures which carry some properties concerning orthogonality types, and he found some interesting characterizations of the Euclidean plane. Moreover, in \cite{martiniantinorms} Martini and Swanepoel briefly discussed angular measures which are proportional to arcs and to areas of sectors of the unit circle.

 Having all this in mind, our purpose is to develop some kind of trigonometry in an arbitrary normed plane, not using any angular measure, but considering a suitable metric extension of the sine function to general normed planes, which was introduced by Szostok \cite{szostok} when studying functional equations. We claim that this generalized sine function has some interesting geometric properties, which we explore throughout this paper.

 The idea of studying trigonometry in normed spaces is not new: functions which somehow play the role of cosine and sine functions in normed spaces were studied by authors like Finsler \cite{finsler}, Busemann \cite{busemann}, and Thompson \cite{thompson} (see Chapter 8 there). The last two references provide definitions for sine functions in terms or area and volume. We show that when dealing with the planar case, our definition based on a metric point of view can be, in some sense, given in terms of Euclidean area (Proposition \ref{prop2}). As observed by Thompson (\cite{thompson}, Section 8.5), Minkowskian trigonometric functions can also be related to solutions of second-order linear differential equations of the type $x'' + f(t)x = 0$ (known as Hill equations).

 In Section \ref{secsine} we define the generalized sine function and present some of its basic properties. Section \ref{secradonsine} is devoted to characterizations of Radon planes using the sine function (in particular, a new proof of D\"{u}velmeyer's characterization of Radon planes via Busemann and Glogovskii angular bisectors presented in \cite{duvelmeyer} is given). In Section \ref{secangle} we prove that if a Radon plane is endowed with an angular measure which provides equal or supplementary measures for angles with the same sine, then this plane has to be Euclidean. In Section \ref{secconstants} sine based constants are defined and used to estimate how  non-Radon or non-Euclidean a normed plane can be (extremal cases are also considered). In Section \ref{secconformal}  we study sine preserving self-mappings, which we call \textit{sine conformal}. Last, a Law of Sines for Radon planes is given in Section \ref{lawofsines}.

 As usual, $(V,||\cdot||)$ denotes a \emph{normed plane}, also called \emph{Minkowski plane}, i.e., an affine plane whose norm $\|\cdot\|$ is determined by a \emph{unit ball} which is an arbitrary
 compact, convex figure centered at the origin; the boundary
 $S = \{x \in V : \|x\| = 1\}$ of this figure is called \emph{unit circle} of $(V, \|\cdot\|)$.
 For such planes we consider three orthogonality types. Given two vectors $x,y \in V$, we say that \\

$\bullet$ $x$ is \textit{Birkhoff orthogonal} to $y$ (denoted by $x \dashv_B y$) when $||x + ty|| \geq ||x||$ for every $t \in \mathbb{R}$,\\

$\bullet$ $x$ is \textit{isosceles orthogonal} to $y$ (denoted by $x \dashv_I y$) when $||x+y|| = ||x-y||$, and \\

$\bullet$ $x$ is \textit{Roberts orthogonal} to $y$ (denoted by $x \dashv_R y$) when $||x + ty|| = ||x - ty||$ for all $t \in \mathbb{R}$. \\

Also, we denote by $[ab]$, $\left.[ab\right>$, and $\left<ab\right>$ the \emph{segment} from $a$ to $b$, the
\emph{half-line} with origin in $a$ and passing through $b$, and the \emph{line} spanned by $a$ and $b$, respectively. 

Good introductions to the geometry of Minkowski planes and spaces are given in the book \cite{thompson}, by the surveys \cite{martini1} and \cite{martini2}, and in the papers \cite{martiniantinorms} and \cite{alonso}.

\section{The generalized sine function} \label{secsine}

When dealing with inner product planes $(V,\left<,\right>)$, the (fairly known) sine between two non-zero vectors $x,y \in V$ is defined by \\
\[ s(x,y) = \sqrt{1 - \frac{\left<x,y\right>^2}{||x||_E^2||y||_E^2}}, \]\\
where the norm comes from the inner product ($||x||_E = \sqrt{\left<x,x\right>}$). Since we are defining the sine function only for non-zero vectors, and since it is clear that $s(\alpha x,y) = s(x,\beta y) = s(x,y)$ for any non-zero real numbers $\alpha$ and $\beta$, we may, for the sake of simplicity, restrict the definition of the sine function to unit vectors. For extending this definition properly to general normed planes, the only suitable tool
is the concept of distances given by the norm. Fortunately, we may describe the sine function of an Euclidean plane easily in terms of distances.

\begin{lemma}\label{lemma1} Let $(V,\left<,\right>)$ be an inner product plane, and denote by $S_E$ the unit circle of the norm derived from the inner product
$\left<,\right>$. If $x, y \in S_E$, then $s(x,y)$ is the distance (in the usual norm) from the origin to the line $l:t\mapsto x + ty$. In other words, $s(x,y) = \inf_{t\in\mathbb{R}}||x+ty||_E$.
\end{lemma}
\noindent\textbf{Proof.} From the standard theory of inner product planes it is clear that this infimum is attained for some $t_0 \in\mathbb{R}$ satisfying $\left<x+t_0y,y\right> = 0$. Hence $t_0 = -\left<x,y\right>$. Now
we have \\
\[ \inf_{t\in\mathbb{R}}||x+ty||_E = ||x+t_0y||_E = \sqrt{\left<x-\left<x,y\right>y,x-\left<x,y\right>y\right>} = \sqrt{1-\left<x,y\right>^2}, \]\\
and this is what we need.
\begin{flushright} $\square$ \end{flushright}
Now we are ready for extending the definition of the sine function in a natural way to general normed planes. \\
\begin{defi}\label{defi1}\normalfont Let $(V,||\cdot||)$ be a normed plane. We define the \textit{sine function} $s:S\times S \rightarrow \mathbb{R}$ by \\
\[s(x,y) =  \inf_{t\in\mathbb{R}}||x+ty||. \]\\
In other words, the sine between $x,y\in S$ is the distance from the origin to the line $l:t\mapsto x +ty$.
\end{defi}

\begin{remark}\label{remark2}\normalfont Despite the fact that we defined the sine function only for unit vectors, we will sometimes abuse of the notation and, for any non-zero $x, y \in V$, denote by $s(x,y)$ the sine between the vectors in the respective directions. In other words, sometimes it is more suitable to think about the sine function as a function defined for directions in $V$.
\end{remark}

The next step is to investigate which properties of the sine function still hold in general normed planes.

\begin{lemma}\label{lemma2} In any normed plane $(V,||\cdot||)$ we have $0 \leq s(x,y) \leq 1$ for every $x,y \in S$. Moreover, $s(x,y) = 0$ if and only if $x = \pm y$, and $s(x,y) = 1$ if and only if $x \dashv_B y$.
\end{lemma}
\noindent\textbf{Proof.} The inequality $||x+ty|| \geq 0$ is immediate, and the other one comes from the equality $||x + ty|| = 1$ for $t = 0$. If $x \neq\pm y$, then the line $l:t\mapsto x+ty$ does not passes through the origin, and so $s(x,y) > 0$. The other direction is obvious. Now we look at the last bi-implication. We have $s(x,y) = 1$ if and only if $||x+ty|| \geq 1 = ||x||$ for every $t \in \mathbb{R}$, i.e., if and only if $x \dashv_B y$. This finishes the proof.
\begin{flushright} $\square$ \end{flushright}
\begin{remark}\label{remark1}\normalfont Notice that what we have done so far clearly holds for normed spaces of dimensions larger than $2$. In general, the sine function is not symmetric. More precisely, in the planar case symmetry of the sine function characterizes Radon planes, and for higher dimensions it characterizes inner product spaces. This will be discussed later.
\end{remark}

\begin{figure}[h]
\centering
\includegraphics{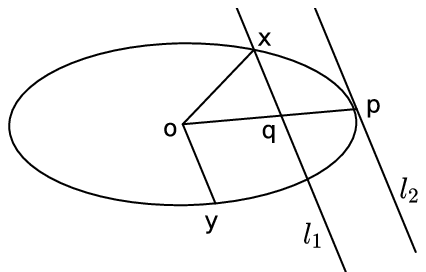}
\caption{$||q|| = s(x,y)$}
\label{fig60}
\end{figure}

 Geometrically we can characterize the sine $s(x,y)$ for linearly independent $x,y \in S$ as follows: the line $l_1:t\mapsto x+ty$ divides the plane into two half-planes. Denote by $H$ the (open) one which does not contain the origin $o$. Hence, a line $l_2$ parallel to $y$ supports $S\cap H$ at a point $p$, say (see Figure \ref{fig60}). Therefore, if the segment $[op]$ intersects $l_1$ at the point $q$, then $||q|| = s(x,y)$.

This follows since distances to lines in direction $y$ are attained by any segments in a direction $z$ with $z \dashv_B y$. In \cite{martini1} it is proved that any normed plane has a \emph{pair of conjugate directions}, i.e., a pair $x,y \in S$ such that $x \dashv_B y$ and $y \dashv_B x$. We will study now some kind of triangle trigonometry for normed planes. Using conjugate diameters, we may also derive some sort of polar coordinates. This is the subject of the next proposition.

\begin{prop}\label{prop3} Let $(V,||\cdot||)$ be a normed plane and assume that $x,y\in S$ are such that $x \dashv_B y$. If $\Delta\mathbf{abc}$ is a triangle such that the segment $[ba]$ is in the direction of $x$ and the segment $[bc]$ is in the direction $y$, then \\
\[ s(c-a,c-b) = \frac{||b-a||}{||c-a||}. \]\\
In particular, if $x,y \in S$ are conjugate diameters, then for any vector $z$, which can be written as $z = \alpha x + \beta y$ for non-negative $a,b \in \mathbb{R}$, we have \\
\[ z = ||z||\left(s(z,y)x + s(z,x)y\right). \]\\
These may be interpreted as polar coordinates for normed planes.
\end{prop}

\begin{figure}[t]
\centering
\includegraphics{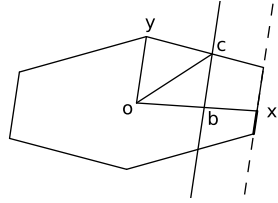}
\caption{$s(c,c-b) = ||b||$}
\label{fig61}
\end{figure}

\noindent\textbf{Proof.} Obviously, we may assume that $a$ is the origin $o$ and that $||a-c|| = 1$. Hence we have to prove that $s(c,c-b) = ||b||$. But this is easy since $s(c,c-b)$ is the distance from the line $\left<cb\right>$ to the origin, and because $x \dashv_B y$, this distance is attained precisely at $b$ (see Figure \ref{fig61}). \\

\begin{figure}[t]
\centering
\includegraphics{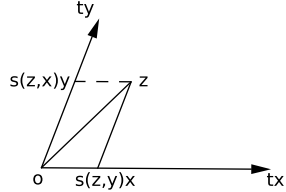}
\caption{Polar coordinates ($z\in S$)}
\label{fig62}
\end{figure}

Assume now that $z = \alpha x + \beta y$ with $\alpha,\beta > 0$, and consider the triangle $\Delta\mathbf{(\alpha x)oz}$ (see Figure \ref{fig62}). By the previous arguments we have $s(z,y) = \frac{||\alpha x||}{||z||} = \frac{\alpha}{||z||}$ and $s(z,x) = \frac{||\beta y||}{||z||} = \frac{\beta}{||z||}$. If $\alpha = 0$ or $\beta = 0$, the result is immediate.
\begin{flushright} $\square$ \end{flushright}

\begin{remark}\label{remark3}\normalfont We can clearly work also with vectors which are not in the quadrant $\mathrm{conv}(\left.[ox\right>\cup\left.[oy\right>)$, by using absolute values. For example, if $z = \alpha x + \beta y$ is a point of $S$ with $\alpha < 0$, then $s(z,x) = |\alpha|$.
\end{remark}

We notice here that triangles which have two sides lying in conjugate directions behave, in the sense of the sine of the corresponding (ordered) angle, like rectangular triangles in the Euclidean plane. Due to this, it would be natural do define a cosine function $c:S\times S \rightarrow\mathbb{R}$ by setting $c(z,x) = s(z,y)$, where $y\in S$ is a conjugate direction of $x$. The problem here is that we can have directions with no conjugate, and if there is a conjugate, it might not be unique. This is the main reason why we work with only one trigonometric function. If we would work with an inner product plane, then we would clearly have $s(z,x)^2 + s(z,y)^2 = 1$ for any orthogonal pair $x,y \in S$ and every $z \in S$. This is not true for an arbitrary normed plane, and thus the ``distortion" in this equality may provide a way to measure how far the plane is from being Euclidean. This is made in Section \ref{secconstants}, but we can already estimate the range of $s(z,x)^2 + s(z,y)^2$.

\begin{figure}[h]
\centering
\includegraphics{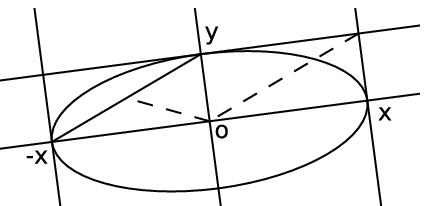}
\caption{Lemma \ref{lemma7}}
\label{fig63}
\end{figure}

\begin{lemma}\label{lemma7} If $x,y \in S$ are conjugate directions, then $\frac{1}{2} \leq s(z,x)^2 + s(z,y)^2 \leq 2$ for all $z \in S$. On the left side, equality holds if and only if $z = \frac{x+y}{2}$ or $z = \frac{y-x}{2}$ (in particular, if one of the segments $[xy]$ or $[y(-x)]$ is contained in the unit circle), and on the right side equality holds if and only if $z = x+y$ or $z = y-x$ (in particular, if $[y(x+y)]$ and $[(x+y)x]$ or $[y(y-x)]$ and $[(y-x)(-x)]$ are contained in $S$).
\end{lemma}
\noindent\textbf{Proof.} It is clear that $s(z,x)^2 + s(z,y)^2 \leq 2$, and that equality holds if and only if $s(z,x) = s(z,y) = 1$. Since a unit vector $z$ of the half-circle $S_1$ from $x$ to $-x$ passing through $y$ can be written as $z = s(z,y)x + s(z,x)y$ or $z = s(z,y)(-x)+s(z,x)y$, it follows that $s(z,x)^2 + s(z,y)^2 = 2$ if and only if $x+y \in S$ or $y-x \in S$. By convexity, if a point $z = \alpha x+\beta y$ belongs to the unit circle, then the intersection of the ray $\left.[oz\right>$ with the segment $[xy]$ is a point $z_0 = \alpha_0x + \beta_0y$ for which $|\alpha_0| \leq |\alpha|$ and $|\beta_0| \leq |\beta|$ hold. Hence, the minimum value for $s(z,x)^2 + s(z,y)^2$ is attained if and only if one of the segments $[xy]$ or $[(-x)y]$ is contained in the unit circle. In this case the value $\frac{1}{2}$ is easily achieved, since $\frac{1}{2} = \min\{\alpha^2+\beta^2:\alpha + \beta = 1\}$. Figure \ref{fig63} illustrates the situation.
\begin{flushright} $\square$ \end{flushright}

Despite the fact that we will discuss the continuity of the sine function later, we will finish this section by using the geometric interpretation to prove that, similarly as in the Euclidean case, the generalized sine function of an arbitrary normed plane attains every value between $0$ and $1$.

\begin{lemma}\label{lemma5} Let $(V,||\cdot||)$ be a normed plane with associated sine function $s:S\times S\rightarrow\mathbb{R}$. Then for every $0 \leq \varepsilon \leq 1$ there exist $x,y \in S$ such that $s(x,y) = \varepsilon$.
\end{lemma}
\noindent\textbf{Proof.} Let $x, y \in S$ be conjugate directions. Let $p \in [ox]$ be such that $||p|| = \varepsilon$, and let $l$ be the line parallel to $y$ and passing through $p$. Hence, if $z$ is any point in $S\cap l$, we have immediately $s(z,y) = \varepsilon$.

\begin{flushright} $\square$ \end{flushright}

\section{Characterizing Radon planes} \label{secradonsine}

Throughout this section we follow the approach given by Martini and Swanepoel in \cite{martiniantinorms}. Thus, $(V,||\cdot||)$ is a Minkowski plane endowed with a non-degenerate symplectic bilinear form $[\cdot,\cdot]:V\times V \rightarrow \mathbb{R}$. We define the \textit{antinorm} of $V$ to be
\[ ||x||_a = \sup_{y\in S}|[x,y]|, \] \\
for every $x\in V$. Moreover, the supremum is attained for some $z \in S$ if and only if $z \dashv_B x$.

\begin{lemma}\label{lemma3} Let $(V,||\cdot||)$ be a normed plane with antinorm $||\cdot||_a$. Then the following statements are equivalent: \\

\noindent\textbf{(a)} \textit{The antinorm is a multiple of the norm}.\\
\noindent\textbf{(b)} \textit{Birkhoff orthogonality is a symmetric relation}.\\
\end{lemma}

For a proof we refer the reader to \cite{martiniantinorms}.
\begin{flushright} $\square$ \end{flushright}

It is well known that a normed plane $(V,||\cdot||)$ which satisfies either \textbf{(a)} or \textbf{(b)} is called a \textit{Radon plane}, see again \cite{martiniantinorms}. \\

In order to characterize Radon planes using only the generalized sine function, we will rewrite it in terms of the antinorm. This rewriting has also some other interesting consequences which show that, in some aspects, 
the sine function behaves for Radon planes as it does in the Euclidean plane.

\begin{prop}\label{prop2} Let $(V,||\cdot||)$ be a normed plane endowed with a nondegenerate symplectic bilinear form $[\cdot,\cdot]$ and with associated antinorm $||\cdot||_a$. Then, for any $x,y \in S$, we have\\
\begin{align}\label{eqsine} s(x,y) = \frac{|[x,y]|}{||y||_a}. \end{align}
\end{prop}
\noindent\textbf{Proof.} If $x = \pm y$, the equality is immediate. If not, then every point of $S\setminus \{y,-y\}$ can be written as $\frac{x+ty}{||x+ty||}$ for some $t \in \mathbb{R}$. Hence\\
\[||y||_a = \sup_{z\in S}|[x,y]| = \sup_{t\in\mathbb{R}}\left|\left[\frac{x+ty}{||x+ty||},y\right]\right| = \sup_{t\in\mathbb{R}}\frac{|[x,y]|}{||x+ty||} = \frac{|[x,y]|}{\inf_{t\in\mathbb{R}}||x+ty||} = \frac{|[x,y]|}{s(x,y)}, \]\\
where one should be aware of the fact that in the second equality the supremum is not attained for $z = \pm y$.
\begin{flushright} $\square$ \end{flushright}

\begin{remark}\label{remark5} \normalfont Notice that if we consider the sine function defined for directions instead of for unit vectors, then formula (\ref{eqsine}) becomes \\
\begin{align}\label{eqsine2} s(x,y) = \frac{|[x,y]|}{||y||_a||x||}. \end{align}
Hence, if we denote by $s_a$ the sine function associated to the antinorm $||\cdot||_a$, we immediately have $s_a(x,y) = s(y,x)$ for any directions $x,y \in V$.
\end{remark}

\begin{coro}\label{coro1} The generalized sine function is continuous in any normed plane.
\end{coro}
This comes straightforwardly from formula (\ref{eqsine}). \\
\begin{flushright} $\square$ \end{flushright}

\begin{prop}\label{prop2} A normed plane $(V,||\cdot||)$ is Radon if and only if its associated sine function is symmetric.
\end{prop}
\noindent\textbf{Proof.} If the norm is Radon, then we may write $||\cdot||_a = \lambda ||\cdot||$ for some $\lambda >0$. Therefore, for any $x,y \in S$ we have \\
\[s(x,y) = \frac{|[x,y]|}{||y||_a} = \frac{|[x,y]|}{\lambda} = \frac{|[y,x]|}{||x||_a} = s(y,x) \,.
\] \\
For the converse we just check whether symmetry of the sine function implies symmetry of Birkhoff orthogonality. In fact, if $x \dashv_B y$, then $s(x,y) = 1$, and thus $s(y,x) = 1$. But this implies $y \dashv_B x$.
\begin{flushright} $\square$ \end{flushright}

\begin{remark}\label{remark6}\normalfont In dimensions $\geq 3$, the symmetry of the generalized sine function implies that the space is an inner product space. In fact, the above shows that the symmetry of the sine function is equivalent to the symmetry of Birkhoff orthogonality, and so the statement follows. For characterizations of inner product spaces via the generalized sine function we refer to \cite{szostok}.
\end{remark}

The next corollary is a kind of Minkowskian analogue of the Euclidean statement ``in a triangle, equal angles yield equal sides".

\begin{coro}\label{coro3} Let $\Delta\mathbf{abc}$ be a triangle in a normed plane $(V,||\cdot||)$. We have $s\left(b-a,c-b\right) = s\left(c-a,c-b\right)$ if and only if $||b-a|| = ||c-a||$.
\end{coro}
\noindent\textbf{Proof.} From the formula (\ref{eqsine2}) we have \\
\[
s\left(b-a,c-b\right) = \frac{|[b-a,c-b]|}{||b-a||.||c-b||_a} \ \mathrm{and} \]\\
\[
s\left(c-a,c-b\right) = \frac{|[c-a,c-b]|}{||c-a||.||c-b||_a}\,.
\] \\
Thus, by $|[b-a,c-b]| = |[c-a,c-b]|$ the desired follows.

\begin{flushright} $\square$ \end{flushright}

The corollary above allows us to characterize isosceles orthogonality (defined in Section \ref{secintro}) via the sine function.

\begin{lemma}\label{lemma8} Let $x,y \in V$ be linearly independent non-zero vectors. Then the following statements are equivalent: \\

\normalfont\noindent\textbf{(i)} $x \dashv_I y$, \\

\noindent\textbf{(ii)} $s(x+y,y) = s(x-y,y)$, \emph{and} \\

\noindent\textbf{(iii)} $s(x+y,x) = s(x-y,x)$.

\end{lemma}
\noindent\textbf{Proof.} By Corollary \ref{coro3}, the proof is immediate.

\begin{flushright} $\square$ \end{flushright}

We continue with deriving some characterizations of Radon planes via ``Euclidean properties'' of the sine function. For the first one we note that we may characterize Radon planes by changing the orientation of the angle in Corollary \ref{coro3}. Of course, this is related to the symmetry of the sine function in these planes.

\begin{lemma}\label{lemma4} In any plane which is not Radon there exists a triangle $\Delta\mathbf{abc}$ for which $||b-a|| = ||c-a||$, but $s\left(c-b,b-a\right) \neq s\left(c-b,c-a\right)$.
\end{lemma}
\noindent\textbf{Proof.} We just have to choose $x, y \in S$ such that $||x||_a \neq ||y||_a$ and consider the triangle $\Delta\mathbf{oxy}$. We have $||x|| = ||y|| = 1$, and it is easy to see that $s\left(y-x,x\right) \neq s\left(x-y,y\right)$.

\begin{flushright} $\square$ \end{flushright}

Let $a,b,c \in V$ be three non-collinear points. The \textit{angle} $\wk\mathbf{abc}$ is the convex hull of the union of the half-lines $\left.[b(a-b)\right>$ and $\left.[b(c-b)\right>$. Having this in mind, we may consider two types of angular bisectors in normed planes:\\

$\bullet$ The \textit{Busemann angular bisector} of $\wk\mathbf{abc}$ is the half-line with origin $b$ in the direction $\frac{a-b}{||a-b||} + \frac{c-b}{||c-b||}$. \\

$\bullet$ The \textit{Glogovskii angular bisector} of $\wk\mathbf{abc}$ is the set of the points $p \in \wk\mathbf{abc}$ equidistant to $\left.[b(a-b)\right>$ and $\left.[b(c-b)\right>$.\\

In \cite{duvelmeyer} D\"{u}velmeyer characterized Radon planes as the only Minkowski planes for which the Glogovskii angular bisector and the Busemann angular bisector coincide for any angle;
see also \cite{martiniantinorms} for a discussion of angular bisectors. We can easily obtain this result by characterizing these angular bisectors in terms of the sine function. This is our next aim.

\begin{prop}\label{prop6} Let $x,y \in S$ be unit vectors which form an angle $\wk\mathbf{xoy}$, and let $z \in S$ be such that the half-line $\left.[oz\right>$ lies in the interior of this angle. Then \\

\noindent\textbf{(a)} $\left.[oz\right>$ is the Glogovskii angular bisector of $\wk\mathbf{xoy}$ if and only if $s(z,x) = s(z,y)$, and \\

\noindent\textbf{(b)} $\left.[oz\right>$ is the Busemann angular bisector of $\wk\mathbf{xoy}$ if and only if $s(x,z) = s(y,z)$.\\

As a consequence we have that the Busemann angular bisector in the norm coincides with the Glogovskii angular bisector in the antinorm (and vice versa). Also, it follows that a normed plane is Radon if and only if these bisectors coincide for any angle.
\end{prop}
\noindent\textbf{Proof.} For \textbf{(a)}, consider a point $p$ within $\wk\mathbf{xoy}$ attaining the distance to the sides of the angle, respectively, at $a \in \left.[ox\right>$ and $b\in\left.[oy\right>$. Then we just have to apply Proposition \ref{prop3} to the triangles $\Delta\mathbf{aop}$ and $\Delta\mathbf{bop}$. \\

We now come to \textbf{(b)}. It is known that the Busemann angular bisector of $\wk\mathbf{xoy}$ is the half-line $\left.[o(x+y)\right>$. Hence the result follows immediately from formula (\ref{eqsine}).\\

The affirmation that Busemann angular bisectors in the norm coincide with Glogovskii angular bisectors in the antinorm (and vice versa) comes now from Remark \ref{remark5}. By Proposition \ref{prop2} it follows that these types of bisectors coincide in any Radon plane. For the converse, assume that the bisectors coincide, and let $x,y \in S$ be in distinct directions. After changing some sign, if necessary, we may consider that the Busemann angular bisector of $\wk\mathbf{xoy}$ is the half-line $\left.[oz\right>$, where $z = \frac{x+y}{||x+y||}$. Since this is also the Glogovskii angular bisector, we have from \textbf{(a)} that $s(z,x) = s(z,y)$, and this gives immediately $||x||_a = ||y||_a$. It follows that the antinorm is a multiple of the norm, and therefore the plane is Radon.

\begin{flushright} $\square$ \end{flushright}

In the Euclidean plane the area of a parallelogram can be calculated as $A = \alpha\beta\sin(\theta)$, where $\alpha$ and $\beta$ are the respective lengths of two consecutive sides and $\theta$ is the angle between these sides. Using a fixed non-degenerate symplectic bilinear form $[\cdot,\cdot]$, the area of the parallelogram with consecutive sides given by $x$ and $y$ is $|[x,y]|$. It turns out that we may characterize Radon planes as being the ones for which this formula still holds.

\begin{prop}\label{prop4} A normed plane $(V,||\cdot||)$ is Radon if and only if there exists a number $\lambda >0$ for which the area of any (ordered) parallelogram $\mathbf{abcd}$ is given by $\lambda||a-b|| \cdot ||a-d||s(v,w)$, where $v$ and $w$ are the unit vectors in the directions $b-a$ and $d-a$, respectively. In this case, rescaling the symplectic bilinear form in such a way that the antinorm coincides with the norm, we will have $\lambda = 1$.
\end{prop}
\noindent\textbf{Proof.} Assume first that $(V,||\cdot||)$ is Radon, and let $\lambda > 0$ be the number such that $||\cdot||_a = \lambda ||\cdot||$. Then, denoting by $A$ the area of the parallelogram $\mathbf{abcd}$, we have \\
\[ A = |[b-a,d-a]| = ||a-b||.||a-d||.|[v,w]| = \lambda||a-b||.||a-d||\frac{|[v,w]|}{||w||_a} = \lambda||a-b||.||a-d||s(v,w). \]\\
For the converse, let $x,y \in S$ be linearly independent unit vectors. Then, considering the parallelogram $\mathbf{oxy(x+y)}$, we have \\
\[ |[x,y]| = \lambda s(x,y) = \lambda\frac{|[x,y]|}{||y||_a} \]\\
for some constant $\lambda >0$. Therefore, $||y||_a = \lambda$. It follows that $||\cdot||_a = \lambda||\cdot||$, and thus the norm is Radon. The remaining part is straightforward.
\begin{flushright} $\square$ \end{flushright}

\section{Angular Measures} \label{secangle}

We follow D\"{u}velmeyer (see \cite{duvelmeyer1}) to introduce angular measures for normed planes in the axiomatic way. The main goal of this section is to prove that the existence of such a measure which is ``coherent'' with the sine function is only possible in the Euclidean plane.

\begin{defi}\label{defi3}\normalfont An \textit{angular measure} in a normed plane $(V,||\cdot||)$ is a Borel measure $\mu$ on the unit circle $S$ satisfying \\

\noindent\textbf{(i)} $\mu(S) = 2\pi$,\\
\noindent\textbf{(ii)} for any Borel set $A \subseteq S$ the equality $\mu(-A) = \mu(A)$ holds,\\
\noindent\textbf{(iii)} for each $v \in S$ we have $\mu(\{v\}) = 0$, and\\
\noindent\textbf{(iv)} any nondegenerate arc of the unit circle has positive measure.\\
\end{defi}

We may define the measure of an angle to be the measure of the arc of $S$ determined by its translate to the origin. In view of this, we may define the measure $\mu(x,y)$  of the angle $\wk\mathbf{xoy}$ between two vectors $x,y \in V\setminus\{0\}$ to be the measure of the smallest arc of the unit circle connecting $\frac{x}{||x||}$ and $\frac{y}{||y||}$. Extending our definition so that opposite half-lines with the same origin also form an angle, it clearly follows from \textbf{(i)}, \textbf{(ii)}, and \textbf{(iv)} that $\mu(x,y) = \pi$ if and only if $\frac{x}{||x||} = -\frac{y}{||y||}$. Standard measure theory gives

\begin{lemma}\label{lemma10} Let $\lambda(S)$ be the length, in the norm, of the unit circle $S$ and consider the arclength parametrization $p:\left[0,\frac{\lambda(S)}{2}\right] \rightarrow S$ of one of the arcs from $x_0$ to $-x_0$, where $x_0 \in S$ is any fixed vector. Then the mapping $t \mapsto \mu(x_0,p(t))$ is continuous.
\end{lemma}

This can be used to derive

\begin{figure}[t]
\centering
\includegraphics{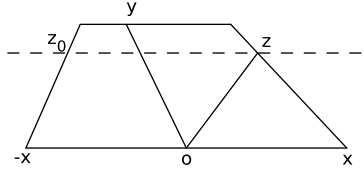}
\caption{$\mu(x,z)=\mu(-x,z_0)$}
\label{fig67}
\end{figure}

\begin{teo}\label{teo3} Let $(V,||\cdot||)$ be a Radon plane. If there exists an angular measure $\mu$ on the unit circle $S$ such that $s(x,y) = s(v,w)$ if and only if $\mu(x,y) = \mu(v,w)$ or $\mu(x,y) + \mu(v,w) = \pi$, then $V$ is the Euclidean plane and $\mu$ is the standard angular measure.
\end{teo}
\noindent\textbf{Proof.} We prove that if such a measure exists, then isosceles orthogonality implies Birkhoff orthogonality. This implication characterizes inner product planes (see \cite{amir}). The first step is to prove that given non-zero vectors $x,y \in V$, we have $\mu(x,y) = \pi/2$ if and only if $x \dashv_B y$. Let us begin by fixing $x,y \in S$ such that $x \dashv_B y$ (and, consequently, $y \dashv_B x$, since we are working with Radon planes). Notice that by the additivity of $\mu$ it follows that $\mu(x,y) + \mu(-x,y) = \pi$. Suppose that $\mu(x,y) > \pi/2$ (the opposite case is analogous). Then, by Lemma \ref{lemma10}, we may choose a point $z$ belonging to the (open) smallest arc from $x$ to $y$ such that $\mu(x,z) = \pi/2$. Assume that the line parallel to $x$ and passing through $z$ intersects the (closed) smallest arc from $y$ to $(-x)$ in a point $z_0$. Thus, the geometric characterization of the sine function given in Section \ref{secsine} shows that $s(-x,z_0) =  s(x,z)$, and hence $\mu(-x,z_0) = \mu(x,z) = \pi/2$. But this is a contradiction, since the union of the smallest arcs from $x$ to $z$ and from $-x$ to $z_0$ is properly contained in an arc joining $x$ to $-x$ (see Figure \ref{fig67}). \\

Now assume that the unit vectors $x$ and $y$ are not Birkhoff orthogonal. Choose $z \in S$ such that $x \dashv_B z$ and assume, changing signs if necessary, that $\wk\mathbf{xoy} \subseteq \wk\mathbf{xoz}$. Since the inclusion is obviously proper, it follows that $\mu(x,y) < \mu(x,z) = \pi/2$. \\

We now prove that if $x,y \in V$ are non-zero vectors which are isosceles orthogonal, then $\mu(x,y) = \pi/2$. By Lemma \ref{lemma8}, for $x \dashv_I y$ we have $s(x+y,x) = s(y-x,x)$ and $s(x+y,y) = s(y-x,y)$.  Notice that $\mu(x,y) + \mu(v,w) = \pi$ if and only if $\mu(x,y) = \mu(-v,w)$. Hence, by the hypothesis and since the vectors $\frac{x+y}{||x+y||}$ and $\frac{y-x}{||y-x||}$ must lie, respectively, in the smallest arcs from $\frac{x}{||x||}$ to $\frac{y}{||y||}$ and from $\frac{y}{||y||}$ to $\frac{-x}{||x||}$, it follows that $\mu(x+y,x) = \mu(y-x,-x)$ and $\mu(x+y,y) = \mu(y-x,x)$. The union of these four respective arcs is an arc which connects $\frac{x}{||x||}$ and $\frac{-x}{||x||}$ (the one which contains $\frac{y}{||y||}$, to be exact). Hence, additivity of $\mu$ gives $\mu(x+y,x) + \mu(x+y,y) = \mu(y-x,y) + \mu(y-x,-x) = \pi/2$. Finally, the equality $\mu(x,y) = \mu(x,x+y) + \mu(x+y,y)$ yields the desired. \\

Therefore, given non-zero vectors $x,y \in V$ such that $x \dashv_I y$, we have $\mu(x,y) = \pi/2$, and this implies $x \dashv_B y$. This shows that $(V,||\cdot||)$ is Euclidean. Moreover, the measure $\mu$ agrees with the Euclidean sine. It follows that $\mu$ is the standard Euclidean measure of angles.

\begin{flushright} $\square$ \end{flushright}

\section{Trigonometric Constants and Extremal Values} \label{secconstants}

This section is devoted to the study of some constants defined via the sine function, and also to characterizations of some types of normed planes for which those constants attain extremal values. We first want to estimate how ``non-Euclidean'' a normed plane $(V,||\cdot||)$ can be. In view of Proposition \ref{prop3} we can do this by measuring how far triangles with two sides in conjugate directions can be from being Pythagorean. For any conjugate directions $x,y \in S$ we define \\
\[ c_E(x,y) := \sup_{z\in S}\left(s(z,x)^2+s(z,y)^2\right) - \inf_{z\in S}\left(s(z,x)^2 + s(z,y)^2\right). \] \\
It is clear that this number is $0$ for any pair of conjugate directions in an inner product plane. The next lemma states that we only need one pair of ``well-behaving'' conjugate directions to characterize such a plane.
\begin{lemma}\label{lemma6} There exists a pair $x,y \in S$ of conjugate directions such that $c_E(x,y) = 0$ if and only if $(V,||\cdot||)$ is an inner product plane.
\end{lemma}
\noindent\textbf{Proof.} If we have conjugate directions $x,y \in S$ such that $c_E(x,y) = 0$, then $s(z,x)^2 + s(z,y)^2 = 1$ for every $z \in S$. Then we simply define an inner product $\left<\cdot,\cdot\right>:V\times V \rightarrow \mathbb{R}$ by setting $\left<x,y\right> =0$ and $\left<x,x\right>=\left<y,y\right>=1$, and we denote by $||\cdot||_E$ the norm given by it. If $z \in S$, we may write $z = \alpha x + \beta y$ for some $\alpha,\beta\in\mathbb{R}$. By Proposition \ref{prop3} we may assume that $|\alpha| = s(z,y)$ and $|\beta| = s(z,x)$. Hence $\alpha^2 + \beta^2 = 1$. Since $||z||_E$ clearly equals $\alpha^2 + \beta^2$, it follows that $||\cdot|| = ||\cdot||_E$, as we wished. The converse is obvious.
\begin{flushright} $\square$ \end{flushright}

 It is natural to ask now what the largest ``distortion" among all pairs of conjugate diameters is. We define the constant
\[ c_E\left(||\cdot||\right) := \sup\{c_E(x,y):x\dashv_B y \ \mathrm{and} \ y\dashv_B x\}. \]\\
For this constant we have the following

\begin{prop}\label{prop5} For any normed plane $(V,||\cdot||)$ we have $0 \leq c_E\left(||\cdot||\right) \leq \frac{3}{2}$. Equality on the left side holds if and only if the norm is derived from an inner product, and equality on the right side holds if and only if the unit circle is an affine regular hexagon.
\end{prop}
\noindent\textbf{Proof.} The left side is immediate due to Lemma \ref{lemma6}. For the right side, notice first that Lemma \ref{lemma7} yields the inequality immediately. By compactness the supremum is in fact a maximum, and hence we have the following: if equality holds, then we have conjugate directions $x,y \in S$ such that $\sup_{z\in S}\left(s(z,x)^2 + s(z,y)^2\right) = 2$ and $\inf_{z\in S}\left(s(z,x)^2+ s(z,y)^2\right) = \frac{1}{2}$. Using again compactness and Lemma \ref{lemma7}, we see that this is only possible if the unit circle is the affine regular hexagon with vertices $\pm x$, $\pm y$ and $\pm (x+y)$, or $\pm (x-y)$.
\begin{flushright} $\square$ \end{flushright}

Recall that Proposition \ref{prop2} states that a Minkowski plane is Radon if and only if its associated sine function is symmetric. Thus, it is natural to estimate how far a norm is from being Radon by calculating the differences between $s(x,y)$ and $s(y,x)$ for pairs $x,y\in S$. We define the respective constant $c_R(||\cdot||)$ by \\
\[ c_R(||\cdot||) = \sup_{x,y\in S}|s(x,y) - s(y,x)|. \]

 In view of Remark \ref{remark5} one can instantly check that $c_R(||\cdot||_a) = c_R(||\cdot||)$. In some sense, according to $c_R$ a norm and its antinorm are ``equally non-Radon" . The next theorem characterizes the rectilinear planes as the ``most non-Radon" planes.

\begin{figure}[h]
\centering
\includegraphics{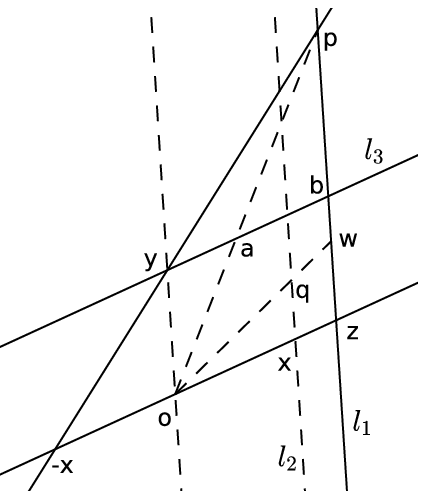}
\label{fig59}
\caption{Proof of Theorem \ref{teo1}}
\end{figure}

\begin{teo}\label{teo1} For any normed plane $(V,||\cdot||)$ we have

\[ c_R(||\cdot||) \leq \frac{1}{2}. \] \\
\noindent Equality holds if and only if the plane is rectilinear.
\end{teo}
\noindent\textbf{Proof.} Let $x,y \in S$ be vectors with $x \neq \pm y$ and assume, without loss of generality, that $s(x,y) \geq s(y,x)$. Denote by $H$ the (open) half-plane determined by the line parallel to $y$ and passing through the origin which contains also $x$, and let $w$ be a point at which a line of direction $y$ supports $S\cap H$. We will denote this supporting line by $l_1$. Then the ray $\left.[ow\right>$ intersects the line $l_2: t \mapsto x + ty$ at the point $q = s(x,y)w$. Let $z$ be the intersection of the ray $\left.[ox\right>$ with the line $l_1$ and assume that $||x-z|| = \beta$ (notice that $\beta = 0$ if and only if $x \dashv_B y$). We may calculate $s(x,y)$ in terms of $\beta$: from the triangles $\Delta\mathbf{owz}$ and $\Delta\mathbf{oqx}$ we have $\frac{||x||}{||x-z||} = \frac{||q||}{||w||}$ (see Figure \ref{fig59}). This gives immediately $s(x,y) = \frac{1}{1+\beta}$.\\

If we denote by $p$ the point where the line $\left<(-x)y\right>$ intersects $l_1$, then it follows that the portion of $S$ from $x$ to $y$ (the one which does not contain $-x$) is contained in $\mathrm{conv}\{y,p,z,o\}$. In particular, this means that the rays of the quadrant $\mathrm{conv}(\left.[oy\right>\cup\left.[ox\right>)$ intersect the unit circle before they intersect $[yp]\cup[pz]$. This follows by convexity and from the fact that $l_1$ supports $S$ (see Figure \ref{teo1}).

We will use this ``location" of that portion of the unit circle to estimate $||y + tx||$ whenever $t \geq 0$. It will be necessary to consider two cases, but first we need to find the intersection $a$ of $\left.[op\right>$ with the line $l_3:t \mapsto y+tx$. If $l_1$ and $l_3$ intersect at $b$ then, by using homothety between triangles, we have \\
\[ \frac{||x+z||}{||y-b||} = \frac{||p-z||}{||p-b||} = \frac{||z||}{||a-b||} \,.
\] \\

 Hence $\frac{2+\beta}{1+\beta} = \frac{1+\beta}{||a-b||}$. Now put $||y-a|| = 1+\beta - ||a-b|| = \frac{1+\beta}{2+\beta}$. Therefore, $a = y + \frac{1+\beta}{2+\beta}x$. Denote by $\lambda_t$ the non-negative number such that $\lambda_t(y+tx)$ is the point where the ray $\left.[o(y+tx)\right>$ intersects the unit circle $S$. Then $||y + tx|| = \lambda_t^{-1}$. If $0 \leq t \leq \frac{1+\beta}{2+\beta}$, a simple calculation shows that the ray $\left[o(y+tx)\right>$ intersects the segment $[yp]$ at the point $\frac{1}{1-t}(y+tx)$. Since the ray intersects the unit circle before it intersects the segment, it follows that $\lambda_t \leq \frac{1}{1-t}$. Hence \\
\[ ||y + tx|| = \lambda_t^{-1} \geq 1-t \geq 1-\frac{1+\beta}{2+\beta} = \frac{1}{2+\beta}. \]\\
If $t \geq \frac{1+\beta}{2+\beta}$, then the ray $\left[o(y+tx)\right>$ intersects the segment $[pz)$ at $\frac{1+\beta}{t}(y+tx)$. Thus $\lambda_t \leq \frac{1+\beta}{t}$, and therefore \\
\[ ||y+tx|| = \lambda_t^{-1} \geq \frac{t}{1+\beta} \geq \frac{1}{2+\beta}. \] \\
We may repeat exactly the same argument to show that this estimate still holds if $t < 0$. Then \\
\[
|s(x,y) - s(y,x)| = s(x,y) - s(y,x)  \leq \frac{1}{1+\beta} - \frac{1}{2+\beta} = \frac{1}{(1+\beta)(2+\beta)} \leq \frac{1}{2} \,,
\] \\
and the desired inequality follows. \\

Notice that if we have equality, then by continuity of the sine function and compactness of $S\times S$ there exists a pair $x, y \in S$ satisfying $s(x,y) - s(y,x) = \frac{1}{2}$. In this case, we have necessarily $\beta = 0$, and therefore $x \dashv_B y$. Also we must have $\lambda_t = 2$, and this happens only when $t = \pm \frac{1}{2}$. If $t = \frac{1}{2}$, then by $\left|\left|y + \frac{1}{2}x\right|\right| = \frac{1}{2}$ it follows that $x + 2y \in S$. Therefore, since $-x$, $y$ and $x+2y$ are collinear points in $S$, we have that $[(-x)(x+2y)]$ is a segment of the unit circle. This, together with central symmetry and the orthogonality $x \dashv_B y$, yields immediately that $S$ is the parallelogram whose vertices are the points $\pm x$ and $\pm (x+2y)$. If $t = -\frac{1}{2}$, we repeat the argument, but will obtain the parallelogram with vertices $\pm x$ and $\pm (-x+2y)$. \\

It remains to prove that the supremum is attained for any rectilinear plane, but this is obvious. If $S$ is the parallelogram with vertices $\pm v$ and $\pm w$, we just have to consider the vectors $v$ and $\frac{v+w}{2}$.
\hfill $\Box$

\begin{flushright} $\square$ \end{flushright}

It is known that the affine regular $(4n+2)$-gons, $n \in \mathbb{N}$, are Radon curves, and that this is never true for the affine regular $(4n)$-gons (this was noticed by Heil in \cite{heil}). We now calculate the distortion $c_R$ for norms whose unit circle is a regular $(4n)$-gon.

\begin{teo}\label{teo2} Let $||\cdot||_{4n}$ denote the norm whose unit circle is given by an affine regular $(4n)$-gon, $n\in\mathbb{N}$. Then we have \\
\[ c_R\left(||\cdot||_{4n}\right) = \left(\sin\frac{\pi}{4n}\right)^ 2. \]
\end{teo}
\noindent\textbf{Proof.} It is clear that the sine function is invariant under linear transformations in the following sense: if $T:V\rightarrow V$ is a linear transformation and $s$ is the sine function associated to a norm whose unit circle is $S$ then, if we denote by $s_T$ the sine function induced by the norm with unit circle $T(S)$, we will have $s(x,y) = s_T(Tx,Ty)$. Consequently, the constant $c_R$ remains the same if we modify the unit circle by an affine transform, and therefore we may use the standard Euclidean $(4n)$-gons to perform the calculations. \\

\begin{figure}[h]
\centering
\includegraphics{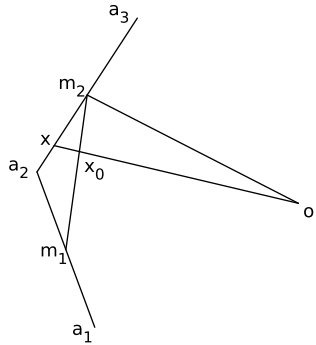}
\label{fig66}
\caption{Estimating $||x||_a$}
\end{figure}

For simplicity, denote by $||\cdot||$ the norm whose unit circle S is a regular $(4n)$-gon. It is known that the anticircle $S_a$ of $S$ is a homothet of the polygon whose vertices are the midpoints of the sides of $S$. In particular, we may rescale $[\cdot,\cdot]$ so that $S_a$ is precisely this polygon. Our first task is to determine, in this case, the minimum and maximum values of $||x||_a$ as $x$ ranges through $S$. First, from $S_a \subseteq S$ it follows that $||x||_a \geq 1$ for every $x \in S$. Since equality holds whenever $x$ is the midpoint of some side of $S$, it follows that $\min_{w\in S}||w||_a = 1$. For the maximum, assume that the ray $\left.[ox\right>$ intersects the segment joining the respective midpoints $m_1$ and $m_2$ of consecutive sides $[a_1a_2]$ and $[a_2a_3]$ of $S$. Since $[m_1m_2] \subseteq S_a,$ we have that this intersection occurs at $x_0 = \frac{x}{||x||_a}$. Thus, the value of $||x||_a$ is the ratio between the Euclidean lengths of the segments $[ox]$ and $[ox_0]$ (see Figure \ref{fig66}).

At this point, our problem becomes a problem of planar Euclidean geometry: we have to find the maximum ratio between lengths of $[ox]$ and $[ox_0]$ as $x$ ranges from $a_2$ to $m_2$ (evidently, we have 
 correspondingly equal values if $x$ varies from $m_1$ to $a_2$). Let $\theta$ be the standard Euclidean measure of the angle $\wk\mathbf{xom_2}$ and assume, without loss of generality, that $||a_2||_E = 1$, where $||\cdot||_E$ is the Euclidean norm. Then basic trigonometry gives \\
\[ \frac{||x||_E}{||x_0||_E} = 1 + \tan\theta\tan\left(\frac{\pi}{4n}\right). \] \\
 Hence the maximum value for $\frac{||x||_E}{||x_0||_E}$ is $1 + \left(\tan\frac{\pi}{4n}\right)^2 = \left(\cos\frac{\pi}{4n}\right)^{-2}$. It follows that $\max_{z\in S}||z||_a = \left(\cos\frac{\pi}{4n}\right)^{-2}$.\\

Now, given $x,y \in S$, we may assume that $s(x,y) \geq s(y,x)$, which is equivalent to $||y||_a \leq ||x||_a$. Due to this we have \\
\[s(x,y) - s(y,x) = s(x,y)\left(1 - \frac{s(y,x)}{s(x,y)}\right) \leq 1 - \frac{||y||_a}{||x||_a} \leq 1 - \frac{\min_{w\in S}||w||_a}{\max_{z\in S}||z||_a} = 1 - \left(\cos\frac{\pi}{4n}\right)^2 = \left(\sin\frac{\pi}{4n}\right)^2 \,.
\]\\
To check that this value is optimal, let $a_1, a_2, ..., a_{4n}$ be the vertices of $S$. Denote by $m_j$ the midpoint of the side $[a_ja_{j+1}]$ (identifying $a_1$ with $a_{4n+1}$). Since $S$ is a regular $(4n)$-gon, the vector $m_n$ is parallel to $[a_{4n}a_1]$, and hence $a_1 \dashv_B m_n$. This gives $|[a_1,m_n]| = ||m_n||_a = 1$ (recall that we rescaled $[\cdot,\cdot]$ in such a way that the unit anticircle is the polygon whose vertices are $m_1, m_2,...m_{4n}$). Moreover, it is clear from the previous calculations that $||a_1||_a = \left(\cos\frac{\pi}{4n}\right)^{-2}$. It follows that $|s(a_1,m_n)-s(m_n,a_1)| = \left(\sin\frac{\pi}{4n}\right)^2$, as desired.
\begin{flushright} $\square$ \end{flushright}

We should note that a geometric constant which can be written in terms of the sine function has been studied previously. It quantifies the difference between isosceles orthogonality and Birkhoff orthogonality and it is defined as: \\
\[ D(X) := \inf\left\{\inf_{\lambda\in \mathbb{R}}||x+\lambda y|| : x,y \in S, \ x \dashv_I y\right\}. \] \\
Clearly, $D(X) = \inf\left\{s(x,y):x,y \in S, \ x\dashv_I y\right\}$. This constant was defined in \cite{wu}. Another reference is \cite{alonso}. 

\section{Sine Conformal Mappings} \label{secconformal}

Conformal mappings are usually defined to be mappings which preserve measures of angles.
For that reason we will call sine preserving self-mappings of $V$ sine conformal. More precisely, we fix this by

\begin{defi}\label{defi2}\normalfont Let $(V,||\cdot||)$ be a normed plane with associated sine function (in the sense of Remark \ref{remark2}) $s: V\setminus\{0\}\times V\setminus\{0\} \rightarrow \mathbb{R}$. A mapping $f:V\rightarrow V$ for which $f(x) \neq 0$ if $x \neq 0$ is said to be \textit{sine conformal} if $s(f(x),f(y)) = s(x,y)$ for every $x,y \in V\setminus\{0\}$.
\end{defi}

Notice that sine conformal mappings, in particular, preserve Birkhoff orthogonality  (see Lemma \ref{lemma2}). It is known that any linear map which preserves Birkhoff orthogonality is a scalar multiple of an isometry (this was established in \cite{blanco}). Using this, we may describe the sine conformal mappings which are linear.

\begin{prop}\label{prop7} A linear map $f:V \rightarrow V$ is sine conformal if and only if it is a scalar multiple of an isometry.
\end{prop}
\noindent\textbf{Proof.} We just have to prove that a scalar multiple of a linear isometry is sine conformal, since the other direction of the bi-implication comes immediately from the previous comments. Let $f = kL$, where $k \neq 0$ and $L: V\rightarrow V$ is a linear isometry. Fix any linearly independent $x,y \in V$ and let $z \in S$ be such that $z \dashv_B y$. Assume that we may write $z = \alpha x + \beta y$. Hence \\
\[
s(x,y) = \frac{|[x,y]|}{||y||_a||x||} =\frac{|[x,y]|}{|[y,z]| \cdot ||x||} = \frac{1}{|\alpha| \cdot ||x||}\,.
\] \\
On the other hand, since $f$ obviously preserves Birkhoff orthogonality, we have that $||f(y)||_a = \left|\left[f(y),\frac{f(z)}{k}\right]\right|$. Hence \\
\[
s(f(x),f(y)) =\frac{|[f(x),f(y)]|}{||f(y)||_a||f(x)||} = \frac{|[f(x),f(y)]|}{|[f(y),f(z)]| \cdot ||x||} = \frac{1}{|\alpha| \cdot ||x||},
\]\\
showing that $f$ is sine conformal.\\
\begin{flushright} $\square$ \end{flushright}

In our next lemma we give a first result concerning the existence of a non-trivial sine conformal mapping.  As a consequence, we will derive a characterization of inner product planes among Radon planes.
A linear map whose eigenvalues are $1$ and $-1$ we call a \textit{reflection}. For the following notice that Roberts orthogonality (defined in Section \ref{secintro}) in the norm is equivalent to Roberts orthogonality in the antinorm. In fact, assume that $x \dashv_R y$. In this case, fixing any $\alpha \in\mathbb{R}$, we have\\
\[ \left|\left[x+\alpha y, \frac{x - ty}{||x - ty||} \right]\right| = \left|\left[x-\alpha y, \frac{x + ty}{||x + ty||} \right]\right| \]\\
for every $t \in \mathbb{R}$. It follows that the continuous maps $z \mapsto |[x+\alpha y,z]|$ and $z \mapsto |[x - \alpha y,z]|$ have the same range when $z$ varies through $S$. Thus, $||x+\alpha y||_a = ||x - \alpha y||_a$.
\begin{lemma}\label{lemma9} Let $f:V \rightarrow V$ be a reflection whose eigenvectors are $x$ and $y$. Then $f$ is sine conformal if and only if we have that $x$ and $y$ are Roberts orthogonal.
\end{lemma}
\noindent\textbf{Proof.} We must have $f(x) = x$ and $f(y) = -y$ or $f(x) = -x$ and $f(y) = y$. In both cases, formula (\ref{eqsine2}) gives that $s(f(v),f(w)) = s(v,w)$ for every $v,w \in V\setminus\{0\}$ if and only if  \\
\[ ||\alpha x + \beta y||_a||\nu x + \xi y|| = ||\alpha x - \beta y||_a||\nu x - \xi y||
\] \\
for any $\alpha,\beta,\nu,\xi \in \mathbb{R}$. Hence $f$ is sine conformal if and only if $x$ and $y$ are Roberts orthogonal in both, the norm and the antinorm. But this is equivalent to $x \dashv_R y$.

\begin{flushright} $\square$ \end{flushright}

\begin{coro}\label{coro4} Let $(V,||\cdot||)$ be a Radon plane. Then $V$ is an inner product plane if and only if for every conjugate pair $\{x,y\} \subseteq S$ the reflection $f:V\rightarrow V$ given by $f(x) = x$ and $f(y) = -y$ is sine conformal.
\end{coro}
\noindent\textbf{Proof.} Radon planes are precisely the planes with the property that for every $x \in S$ there exists a $y \in S$ such that $x \dashv_B y$ and $y \dashv_B x$. Hence the hypothesis that the reflection given by $f(x) = x$ and $f(y) = -y$ is sine conformal whenever $x$ and $y$ are conjugate gives that for every $x \in S$ there exists some $y \in S$ for which $x \dashv_R y$. This characterizes inner product planes (see \cite{alonso}). The converse is obvious.

\begin{flushright} $\square$ \end{flushright}

\section{The Law of Sines} \label{lawofsines}

In this final part we present a few results related to the Law of Sines in Radon planes.

\begin{teo}[Law of Sines]\label{teolaw} Let $(V,||\cdot||)$ be a Radon plane, and let $x,y,z \in S$ be non-collinear. Then in the triangle $\Delta\mathbf{xyz}$ we have\\
\begin{align}\label{sinelaw} \frac{||x-y||}{s(x-z,y-z)} = \frac{||y-z||}{s(x-y,x-z)} = \frac{||x-z||}{s(y-z,x-y)}. \end{align}\\
If this ratio equals $2$ for any triple of distinct unit vectors, then $V$ is an inner product plane. \\
\end{teo}
\noindent\textbf{Proof.} The equality follows immediately from (3.2). Assume now that the ratio is $2$ for any $x,y,z \in S$. If $x,y \in S$, then, looking to the triangle $\Delta\mathbf{xy(-x)}$, we have \\
\[ \frac{||2x||}{s(x+y,x-y)} = 2, \]\\
and hence $s(x+y,x-y) = 1$. It follows that $(x+y) \dashv_B (x-y)$ whenever $x,y \in S$. This property characterizes inner product planes (see \cite{alonso}, the comment below Theorem 4.20).

\begin{flushright} $\square$ \end{flushright}

\begin{remark}\normalfont Notice that to characterize inner product planes it is enough to demand only that the ratio (\ref{sinelaw}) equals $2$ for any triangle inscribed in the unit circle having a diameter as one of its sides. In other words, this means that ``any inscribed angle which is opened to a diameter is a right angle". For a related classification of triangles in normed planes (having then obtuse, right, and acute ones) we refer to the paper \cite{AMS}.
\end{remark}

The question that arises is whether or not for any Radon plane we have a triangle inscribed in the unit circle for which the ratio expressed in Theorem \ref{teolaw} equals $2$. To answer this question we enunciate (a part of) the main theorem of \cite{benitez}.

\begin{teo}\label{teo4} Let $(V,||\cdot||)$ be a normed space. For any $x,y \in V\setminus\{0\}$ there exists a unique number $\alpha = \alpha(x,y) > 0$ such that $(x+\alpha y) \dashv_B (x - \alpha y)$. Moreover, the function which associates each pair $(x,y) \in V\setminus\{0\}\times V\setminus\{0\}$ to the number $\alpha(x,y) \in \mathbb{R}$ is continuous.
\end{teo}
\begin{flushright} $\square$ \end{flushright}

\begin{prop}\label{prop8} In any Radon plane, there exist vectors $x,y \in S$ with different directions such that $(x+y) \dashv_B (x-y)$.
\end{prop}
\noindent\textbf{Proof.} Fix arbitrary $v,w \in S$ with $v \neq \pm w$ and let $\alpha:V\setminus\{0\}\times V\setminus\{0\} \rightarrow \mathbb{R}$ be as in Theorem \ref{teo4}. Define $f:[0,1] \rightarrow \mathbb{R}$ by \\
\[ f(\lambda) = \alpha\left(\frac{(1-\lambda)v+\lambda w}{||(1-\lambda)v + \lambda w||}, \frac{(1-\lambda)w+\lambda (-v)}{||(1-\lambda)w + \lambda (-v)||}\right). \]\\
Notice that $f$ is continuous. Moreover, since the plane is Radon, we have $f(0) = \alpha(v,w) = \alpha(w,-v)^{-1} = f(1)^{-1}$. It follows from the Intermediate Value Theorem that there exists some $\lambda_0 \in [0,1]$ such that $f(\lambda_0) = 1$. Setting $x = \frac{(1-\lambda_0)v + \lambda_0w}{||(1-\lambda_0)v + \lambda_0w||}$ and $y = \frac{(1-\lambda_0)w + \lambda_0(-v)}{||(1-\lambda_0)w + \lambda_0(-v)||}$, we have the desired. \\

\begin{flushright} $\square$ \end{flushright}

\begin{coro} In any Radon plane there exists a triangle inscribed in the unit circle for which the ratio (\ref{sinelaw}) attains the ``Euclidean value'' $2$.
\end{coro}
\noindent\textbf{Proof.} Let $x,y \in S$ be as in Proposition \ref{prop8} and consider the triangle $\Delta\mathbf{xy(-x)}$. The statement follows immediately.

\begin{flushright} $\square$ \end{flushright}

\begin{remark}\label{remark7}\normalfont Inspired by formula (\ref{eqsine2}), one can easily verify a weaker Law of Sines for non-Radon planes. Given a triangle $\Delta\mathbf{abc}$ in such a plane, we have the equality \\
\[ \frac{||c-a||}{s(b-a,c-b)} = \frac{||b-a||}{s(c-a,c-b)},  \] \\
and the two analogous equalities (obtained by interchanging suitably the vertices). Only in a Radon plane the three equalities yield the same value for any given triangle. Moreover, notice that Proposition \ref{prop3} and Corollary \ref{coro3} can be seen as special cases of this law. 
\end{remark}

\end{document}